\documentclass[12pt]{article}       % onecolumn (second format)
\usepackage{graphicx}
 \usepackage{mathptmx}      % use Times fonts if available on your TeX system
 \usepackage{amsmath, amssymb, amsthm}
\newtheorem{thm}{Theorem}[section]
\newtheorem{ass}[thm]{Assumption}
\newtheorem{lem}[thm]{Lemma}

\newtheorem{corollary}[thm]{Corollary}
\newtheorem{defn}[thm]{Definition}

\begin{document}

\title{A consistent BGK model with velocity-dependent collision frequency for gas mixtures%\thanks{Grants or other notes
%about the article that should go on the front page should be
%placed here. General acknowledgments should be placed at the end of the article.}
}
%\subtitle{Do you have a subtitle?\\ If so, write it here}

%\titlerunning{Short form of title}        % if too long for running head

\author{J.~ Haack\thanks{Los Alamos National Laboratory, Los Alamos, NM 87545, USA, haack@lanl.gov } \and C.~Hauck\thanks{Oak Ridge National Laboratory, 1 Bethel Valley Road, Bldg. 5700,
Oak Ridge, TN 37831-6164,
USA,
 hauckc@ornl.gov} \and C.~Klingenberg\thanks{klingen@mathematik.uni-wuerzburg.de,
 Universität Würzburg,
Emil-Fischer-Str. 40,
97074 Würzburg, Germany} \and M.~Pirner\thanks{marlies.pirner@mathematik.uni-wuerzburg.de,
Universität Würzburg,
Emil-Fischer-Str. 40,
97074 Würzburg, Germany} \and S.~Warnecke\thanks{sandra.warnecke@mathematik.uni-wuerzburg.de, Universität Würzburg,
Emil-Fischer-Str. 40,
97074 Würzburg, Germany} %etc.
}

%\authorrunning{Short form of author list} % if too long for running head

%\footnote{Jeffrey Haack \\
%                     %  \\
%%             \emph{Present address:} of F. Author  %  if needed
%           \\ \\ Cory Hauck
%            \\
%  \\ \\ Christian Klingenberg, Marlies Pirner, Sandra Warnecke \\ klingen@mathematik.uni-wuerzburg.de \\ marlies.pirner@mathematik.uni-wuerzburg.de \\ 
%%Gebäude: 40 (Mathematik Ost)
%%Raum: 03.011
%%Telefon: +49 931 31-85045
%%Fax: +49 931 31-80944 
%\\ %\and Marlies Pirner, Universität Wien, Oskar-Morgenstern-Platz 1, 1090 Wien \\marlies.pirner@univie.ac.at
%% \and Sandra Warnecke \at Emil-Fischer-Str. 40 , 97074 Würzburg, Germany \\  sandra.warnecke@mathematik.uni-wuerzburg.de 
%}

%\date{Received: date / Accepted: date}
% The correct dates will be entered by the editor

\maketitle

\begin{abstract}
We derive a multi-species BGK model with velocity-dependent collision frequency for a non-reactive, multi-component gas mixture. The model is derived by minimizing a weighted entropy under the constraint that the number of particles of each species, total momentum, and total energy are conserved. We prove that this minimization problem admits a unique solution for very general collision frequencies. Moreover, we prove that the model satisfies an H-Theorem and characterize the form of equilibrium. \\ \\
{\bf Keywords: } multi-fluid mixture, kinetic model, BGK approximation \and plasma physics, velocity-dependent collision frequency, entropy minimization
% \PACS{PACS code1 \and PACS code2 \and more}
% \subclass{MSC code1 \and MSC code2 \and more}
\end{abstract}

\section{Introduction}

In this paper, we present a BGK-type model for gas mixtures that, in the case of two species, takes the form
 \begin{align} \begin{split} \label{BGK2}
\partial_t f_1 + v \cdot \nabla_x  f_1   &= \nu_{11}  (M_{11} - f_1) + \nu_{12}  (M_{12}- f_1),
\\ 
\partial_t f_2 + v \cdot \nabla_x  f_2 &=\nu_{22}  (M_{22} - f_2) + \nu_{21}  (M_{21}- f_2), 
\end{split}
\end{align}
along with appropriate boundary and initial conditions.
Here $f_1 = f_1(x,v,t)$ and $f_2 = f_2(x,v,t)$ are the number densities of species of mass $m_1$ and $m_2$, respectively, with respect to the phase space measure $dxdv$; $x\in \mathbb{R}^3$ is the position coordinate of phase space; $v\in \mathbb{R}^3$ is the velocity coordinate; and $t\geq 0$ is time.  The relaxation operator on the right hand side of \eqref{BGK2} involves target functions of the form
\begin{align}
M_{kj} =  \exp(m_k\lambda^{kj}_0 + m_k \lambda^{kj}_1  \cdot  v + m_k\lambda^{kj}_2   |v|^2),
\label{form1}
\end{align}
which depend on parameters $\lambda^{kj}=(\lambda^{kj}_0, \lambda^{kj}_1, \lambda^{kj}_2) \in \mathbb{R}\times \mathbb{R}^3 \times \mathbb{R}^+$, and (non-negative) collision frequencies $\nu_{kj}$.  These parameters depend implicitly on $f_1$ and $f_2$, and once specified, determine the BGK operator.

The purpose of the relaxation operator in \eqref{BGK2} is to provide an approximation of the multi-species Boltzmann collision operator that is more computationally tractable, but still maintains important structural properties.  In the single-species case, the original BGK model \cite{BGK} serves this purpose.  In particular, it has the same collision invariants as the Boltzmann operator (which lead to conservation of number, momentum, and energy) and it satisfies an H-Theorem.  In the multi-species case, these requirements are not as straight-forward to satisfy, but it can be done.  There are many BGK models for gas mixtures proposed in the literature \cite{gross_krook1956,hamel1965,Greene,Garzo1989,Sofonea2001,Pirner,haack,Bobylev,AndriesAokiPerthame2002}, many of which satisfy these basic requirements and, in addition, are able to match some prescribed relaxation rates and/or transport coefficients that come from more complicated physics models or from experiment.  Many of these approaches have been extended to accommodate ellipsoid statistical (ES-BGK) models, polyatomic molecules, chemical reactions or quantum gases; see for example \cite{Pirner2,Todorova,Groppi,Pirner5,Pirner6,Bisi,Bisi2,Pirner9}.

A common feature of all the models mentioned above is that
they only allow for collision frequencies which are independent of the microscopic velocity $v$ of the particles \cite{StruchtrupBook}. 
However, the collision frequencies in principle should depend on the microscopic velocity, which is typically neglected for the reason of simplicity. In the case of neutral gases, velocity independent collision frequency leads to transport properties in the fluid regime that are inconsistent with the full kinetic collision operator, e.g., the Prandtl number. Models such as the ES-BGK model and the Shakov model make changes to the target Maxwellian to provide  extra degrees of freedom to the system, but still retain the constant collision frequency assumption. Some attempts have been proposed to re-introduce velocity dependence in the case of variable hard spheres interactions for neutral gases \cite{MieussensStruchtrup}, for which velocity-dependent collision frequencies are monotonically increasing and are well-defined. For particles interacting with long-ranged Coulomb interactions, i.e., a plasma, the canonical collision rate definition using the cross section is no longer well defined due to a singularity at a zero relative velocity. A velocity-dependent collision frequency is instead defined by the momentum transfer cross section without an integral, which results in a collision frequency that is decreasing in the limit of large relative velocities \cite{LeeMore,KrallTrivelpiece}.

In this paper, we derive a model of the form \eqref{BGK2} that allows for  velocity-dependent collision frequencies.  Our derivation includes as a by-product the single-species BGK model with velocity-dependent collision frequency that was proposed in \cite{Struchtrup}.  We identify target functions that are consistent with the conservation laws for \eqref{BGK2} and satisfy an entropy minimization principle.   In particular, \textit{intra-species} collisions (between the same species) should preserve mass, momentum, and energy within a species; that is, 
 \begin{align}
\int m_k \nu_{kk} \begin{pmatrix}
1 \\  v \\  |v|^2
\end{pmatrix} (M_{kk} - f_k) dv  = 0, \quad k \in \{1,2\}.
\label{constraints_intra}
\end{align}
Meanwhile \textit{inter-species} collisions (between different species) should preserve the mass of each species, but only the combined momentum and energy of both; that is, 
\begin{align}
\begin{split}
\int m_1 \nu_{12}  (M_{12} - f_1) dv = 0, &\quad \int m_2 \nu_{21}  (M_{21} - f_2) dv = 0 \\
  \int m_1 \nu_{12} 
  \begin{pmatrix}
 v \\  |v|^2
\end{pmatrix} (M_{12} - f_1) dv &+ 
\int m_2 \nu_{21} 
 \begin{pmatrix}
 v \\  |v|^2
\end{pmatrix} 
(M_{21} - f_2)dv = 0.
\label{constraints_inter}
\end{split}
\end{align}
When the collision frequencies are independent of $v$, the integrals in \eqref{constraints_intra} and \eqref{constraints_inter} can be computed explicitly, thereby providing relationships between the parameters $\lambda^{kj}$ and the moments of $f_1$ and $f_2$ with respect to $\{1,v,|v|^2\}$.  In the single-species case, this relationship defines the target function as the Maxwellian associated to $f$, while in the multi-species case, additional constraints must be imposed.  However, when the collision frequencies depend on $v$, the aforementioned integrals are not always computable in closed form and the relationship between the parameters $\lambda^{kj}$ and the moments of $f_1$ and $f_2$ with respect to $\{1,v,|v|^2\}$ cannot be written down analytically. 

In spite of the difficulty of relating the target parameters to the moments of the kinetic distributions, the entropy minimization formulation can be still used to establish a unique set of parameters, under the conditions $\lambda_1^{12}= \lambda_1^{21}$ and $\lambda_2^{12} = \lambda_2^{21}$.  We do so by adapting the strategy from \cite{Junk} to fit the current setting.  While a more abstract approach based solely on convex optimization tools can also be used \cite{Borwein}, we follow \cite{Junk} because it provides a more concrete connection to the application at hand.  Our proof provides a rigorous justification for the target function used in \cite{Struchtrup} for the single species case.  It also leads to an H-Theorem for the multi-species system \eqref{BGK2}.

The remainder of the paper is organized as follows. %\cdh{Will need to revisit this paragraph after finishing the paper.} 
In Section \ref{sec2}, we motivate the choice of the target Maxwellians as solutions of minimization problems of the entropy under certain constraints. In Section \ref{sec3}, we prove existence and uniqueness of the minimization problems. In Section \ref{sec4}, we prove consistency of the model meaning that it satisfies the conservation properties, the H-Theorem and Maxwell distributions with equal mean velocity and temperature in equilibrium. In Section \ref{sec5}, we briefly summarize the straightforward extension to the case of $N$ species, still with binary interactions.

\section{The structure of the target functions}

\label{sec2}
In this section, we motivate the form of the target functions in \eqref{form1}.  
It will be convenient in what follows to define the strictly convex function
\begin{align}
\label{eq:entropy}
h(z)= z \ln z -z, \quad z>0,
\end{align}
and the vector-valued function 
\begin{align}
a^k(v) = \begin{pmatrix}
 a_0^k (v) \\  a_1^k (v) \\  a_2^k (v) 
\end{pmatrix}
= \begin{pmatrix}
 m_k \\  m_k v \\  m_k |v|^2 
\end{pmatrix}.
\end{align}
Since $h$ is convex and $h'(z) = \ln(z)$, it follows that
\begin{equation}
\label{eq:h_tangent}
h(x) \geq h(y) + \ln(y)(x-y), \quad \forall ~ y,x \in \mathbb{R}^+.
\end{equation}

\subsection{The one species target Maxwellians}
We seek a solution of the weighted entropy minimization problem
\begin{align}
\min_{g \in \chi_k} \int \nu_{kk} h(g) dv, \quad k \in\{1,2\},
\label{min}
\end{align}
where 
\begin{align}
\chi_k= \left\lbrace g ~ \Big |~ g \geq 0,\, \nu_{kk} (1+ |v|^2) g \in L^1(\mathbb{R}^3),\, \int \nu_{kk} a^k(v)
(g - f_k) dv = 0 \right\rbrace.
\end{align}
The choice of the set $\chi_k$ ensures the conservation properties \eqref{constraints_intra} for intra-species collisions.  The motivation for weighting the usual objective by the collision frequencies in \eqref{min} is that the ansatz will take the form \eqref{form1}. Indeed, by standard optimization theory, any critical point $(M_{kk}, \lambda^{kk})$ of the Lagrange functional $L_k \colon \chi_k \times \mathbb{R}^5 \to \mathbb{R}$, given by
\begin{align}
L_k(g, \alpha) = \int \nu_{kk} h(g) dv - \alpha \cdot \int \nu_{kk} a^k(v)
(g - f_k) dv , 
\end{align}
satisfies the first-order optimality condition
\begin{align}
\frac{\delta L_k}{\delta g} (M_{kk}, \lambda^{kk}) =  \nu_{kk} ( \ln M_{kk} - \lambda^{kk} \cdot a^k(v) ) = 0,
\end{align}
which implies then that
\begin{align}
M_{kk} = \exp{( \lambda^{kk} \cdot a^k)} =  \exp{( m_k\lambda^{kk}_0  +  m_k \lambda^{kk}_1 \cdot  v + m_k\lambda^{kk}_2   |v|^2)}.
\label{form3}
\end{align} 
In Section \ref{sec:single_species_rigorou}, we prove in a rigorous way that there exists a unique function of the form \eqref{form3} that satisfies these constraints. 

\begin{thm}
\label{thm:unique_min}
Suppose that there exists $\lambda^{kk} \in \mathbb{R}\times \mathbb{R}^3 \times \mathbb{R}$ such that the function $M_{kk}$ given in \eqref{form3} is an element of $\chi_k$.  Then $M_{kk}$ is the unique minimizer of \eqref{min}.
\end{thm}
\begin{proof}
According to \eqref{eq:h_tangent}
\begin{equation}
h(g) 
\geq  h(M_{kk}) + \lambda^{kk} \cdot a^k (g-M_{kk}),
\end{equation}
point-wise in $v$. Thus, because $\nu_{kk} \geq 0$, it follows that for all $g \in \chi_k$,
\begin{align}
\int \nu_{kk} h(g) dv 
\geq \int \nu_{kk} h(M_{kk}) dv + \int \nu_{kk}\lambda^{kk} \cdot a^k (g-M_{kk}) dv
= \int \nu_{kk} h(M_{kk}) dv
\end{align}
Hence $M_{kk}$ is a minimizer of \eqref{min}, and uniqueness follows directly from the strict convexity of $h$.
\end{proof}

\subsection{The mixture target Maxwellians}
\label{subsec:mixture_formal_derivation}

For interactions between species, we seek a solution of the weighted entropy minimization problem
\begin{align}
\label{eq:entropy_min_cross_terms}
\min_{g_{1}, g_{2} \in \chi_{12}}  \int \nu_{12} h(g_{1}) dv + \int \nu_{21} h(g_{2}) dv , 
\end{align}
where
\begin{align}
\label{eq:chi12_def}
\begin{split}
\chi_{12}= \Bigg\lbrace (g_{1}, g_{2})~\Big|& ~ g_{1}, g_{2} >0,\, 
\nu_{12} (1+ |v|^2) g_{1},\,
\nu_{21} (1+ |v|^2) g_{2}  \in L^1(\mathbb{R}^3),\\ 
&\int m_1 \nu_{12} g_{1} dv = \int m_1 \nu_{12} f_1 dv, \quad 
\int  m_2 \nu_{21} g_{2} dv = \int m_2 \nu_{21} f_2 dv, \\
&\int m_1 \nu_{12} 
\begin{pmatrix}
 v \\  |v|^2
\end{pmatrix} (g_{1} - f_1) dv  + 
\int m_2 \nu_{21} \begin{pmatrix}
 v \\ |v|^2
\end{pmatrix} (g_{2} - f_2) dv = 0 \Bigg\rbrace.
\end{split}
\end{align}
Here, $\chi_{12}$ is chosen such that the constraints \eqref{constraints_intra} for inter-species collisions are satisfied.  
Similar to the case of intra-species collisions, we consider the Lagrange functional $L \colon \chi \times \mathbb{R} \times \mathbb{R} \times \mathbb{R}^3 \times \mathbb{R} \to \mathbb{R} $
\begin{align}
\begin{split}
L(g_1, g_2, \alpha_0^1,\alpha_0^2, \alpha_1, \alpha_2) &= \int \nu_{12} h(g_{1}) dv + \int \nu_{21} h(g_{2}) dv \\
&-\alpha_0^1 \int m_1 \nu_{12} (g_{1} - f_1 )dv
- \alpha_0^2 \int m_2 \nu_{21} (g_{2} - f_2 )dv \\ 
&- \alpha_1 \cdot \left( \int m_1 \nu_{12}  v (g_1 - f_1) dv + \int m_2 \nu_{21}  v (g_2 - f_2) dv \right)\\
&- \alpha_2 \left( \int m_1 \nu_{12}  |v|^2 (g_1 - f_1) dv +  \int m_2 \nu_{21}  |v|^2 (g_2 - f_2) dv \right).\\
\end{split}
\end{align}%
Any critical point $(M_{12}, M_{21}, \lambda_0^1,\lambda_0^2, \lambda_1, \lambda_2)$ of $L$ satisfies the first-order optimality conditions 
\begin{align}
\frac{\delta L}{\delta g_1}(M_{12}, M_{21}, \lambda_0^1,\lambda_0^2, \lambda_1, \lambda_2)
 = \nu_{12} (  \ln M_{12} - \lambda^{12} \cdot a^1(v) )  =0, \\
 \frac{\delta L}{\delta g_2}(M_{12}, M_{21}, \lambda_0^1,\lambda_0^2, \lambda_1, \lambda_2)
 = \nu_{21} (  \ln M_{21} - \lambda^{21} \cdot a^2(v) )  =0,
\end{align}
where $\lambda^{12} = (\lambda_0^1, \lambda_1, \lambda_2)$ and $\lambda^{21} = (\lambda_0^2, \lambda_1, \lambda_2)$.  Therefore
\begin{align}
M_{12} &= \exp(\lambda^{12} \cdot a^1(v)) = \exp \left(m_1 \lambda_0^{12} + m_1  \lambda_1 \cdot v  + m_1 \lambda_2  |v|^2 \right) \\
M_{21} &= \exp(\lambda^{21} \cdot a^2(v)) = \exp \left(m_2 \lambda_0^{21} + m_2  \lambda_1 \cdot v  + m_2  \lambda_2  |v|^2 \right).
\end{align}
%We define $\lambda_1^{12} = m_1 \lambda_1 + \frac{\bar{\lambda_1}}{\bar{n}_1 \bar{n}_2}, \lambda_1^{21} = m_2 \lambda_1 - \frac{\bar{\lambda_1}}{\bar{n}_1 \bar\lambda_2^{21} = -\frac{3 m_2}{\bar{n}_1 \bar{n}_2} \tilde{\lambda}_4 + \frac{m_2}{2} \lambda_2$. 
Since we only require conservation of the {\em combined} momentum and kinetic energy, there is only one Lagrange multiplier for the momentum constraint and one Lagrange multiplier for the energy constraint. Therefore, $\lambda_1^{12} = \lambda_1^{21}$ and $\lambda_2^{12} = \lambda_2^{21}$ in \eqref{form1}.  When the collision frequency is constant, this restriction is the same as the one used in \cite{haack}, but more restrictive than the model in \cite{Pirner}.

In the next section, we prove the existence of functions of the form \eqref{form1} that satisfy the constraints in \eqref{constraints_intra} and \eqref{constraints_inter}. As in the single species case, it follows that these functions are unique minimizer of the corresponding minimization problem.

\begin{thm}
\label{thm2.2}
	Assume that there exist 
	$\lambda_0^{12} \in \mathbb{R}$, 
	$\lambda_0^{21} \in \mathbb{R}$, 
	$\lambda_1^{12} = \lambda_1^{21} \in  \mathbb{R}^3$, and 
	$\lambda_2^{12} = \lambda_2^{21} \in  \mathbb{R}$
	such that the pair $(M_{12},M_{21})$, where $M_{kj}$ is defined in \eqref{form1}, is an element of $\chi_{12}$.  Then $(M_{12},M_{21})$ is the unique minimizer of \eqref{eq:entropy_min_cross_terms}.
\end{thm}

\begin{proof}
	According to \eqref{eq:h_tangent} 
	\begin{equation}
	h(g) 
	\geq  h(M_{kj}) + \lambda^{kj} \cdot a^k (g-M_{kj}),
	\end{equation}
	point-wise in $v$, for any measurable function $g$ and $k,j \in \{1,2\}$. 
	Therefore, since $\nu_{kj} \geq 0$, it follows that for any measureable functions $g_1$ and $g_2$,
	\begin{align}
	\label{eq:convex_h1h2}
	\int \nu_{12}h(g_1) dv + \int \nu_{21}& h(g_2)  dv
		\geq  \int \nu_{12} h(M_{12}) dv 
			+ \int \nu_{21} h(M_{21})dv  \nonumber \\
			& \quad + \lambda^{12} \cdot\int \nu_{12}  a^1 (g_1-M_{12})dv
			+ \lambda^{21} \cdot \int \nu_{21}  a^2 (g_2-M_{21})dv.
	\end{align}
	Since $\lambda_1^{12}=\lambda_1^{21}$ and $\lambda_2^{12}=\lambda_2^{21}$,
	\begin{align}
\begin{split}
	\lambda^{12} \cdot\int \nu_{12} a^1& (g_1-M_{12}) dv
		+ \lambda^{21} \cdot \int \nu_{21}  a^2 (g_2-M_{21})dv \\
	&= \lambda_0^{12} \int \nu_{12} m_1  (g_1-M_{12})dv
		+ \lambda_0^{21} \int \nu_{21} m_2 (g_2-M_{21})dv \\
		& \quad + \lambda_1^{12} \cdot \left(\int \nu_{12} m_1 v (g_1-M_{12}) dv
		  + \int \nu_{21} m_2 v (g_2-M_{21}) dv \right) \\
		  & \quad + \lambda_2^{12} \left(\int \nu_{12}  m_1 |v|^2 (g_1-M_{12}) dv
		  + \int \nu_{21}  m_2 |v|^2 (g_2-M_{21}) dv \right) .
\end{split}
	\end{align}
	If $(g_1,g_2)$ and $(M_{12},M_{21})$ are elements of  $\chi_{12}$, then the constraints in \eqref{eq:chi12_def} imply that each of the terms above is zero.  
	In such cases, \eqref{eq:convex_h1h2} reduces
	\begin{align}
	\int \nu_{12}h(g_1) dv + \int \nu_{21} h(g_2)  dv
	\geq  \int \nu_{12} h(M_{12}) dv 
	+ \int \nu_{21} h(M_{21})dv,
	\end{align}
	which shows that $(M_{12},M_{21})$ solves \eqref{eq:entropy_min_cross_terms}.  Since the collision frequencies $\nu_{12}$ and $\nu_{21}$ are non-negative and $h$ is strictly convex, it follows that this solution is unique.
	
\end{proof}

\section{Existence and uniqueness of the target Maxwellians}
\label{sec3}
In this section, we prove the existence of the multipliers $\lambda^{11}$, $\lambda^{22}$, $\lambda^{12}$ and  $\lambda^{21}$ such that the single-species targets $M_{11}$ and $ M_{22}$ satisfy \eqref{constraints_intra} and the mixture targets $M_{12}$ and $M_{21}$ satisfy \eqref{constraints_inter}.  We follow closey the strategy laid out in \cite{Junk}, although some variations will be needed to account for the velocity-dependent collision frequencies and the mixture targets.  

Throughout the paper, we denote a distribution function of exponential form by
\begin{align}
\exp^k_{\lambda}(v):= \exp (\lambda \cdot a^k(v)), \quad \lambda=(\lambda_0, \lambda_1, \lambda_2) \in \mathbb{R}^5.  
\label{exp}
\end{align}
and let
\begin{align}
D_{kj}= \lbrace g \geq 0 \mid \nu_{kj}(1+ |v|^2) g \in L^1(\mathbb{R}^3), \, g \not\equiv 0 \rbrace, \quad
\Lambda^{kj} = \lbrace \lambda \in \mathbb{R}^5 \mid \exp^k_{\lambda} \in D_{kj} \rbrace .
\end{align}
For any $g \in D_{kj}$ the moment map $\mu^{kj}$ is given by
\begin{equation}
    \mu^{kj}(g) = \begin{pmatrix}
     \mu^{kj}_0 \\ \mu^{kj}_1 \\ \mu^{kj}_2 
    \end{pmatrix}(g)
    = \int \nu_{kj} a^k(v) g(v) dv .
\end{equation}
We make the following assumptions about the collision frequencies.
\begin{ass}
\label{ass:nu}
Each frequency $\nu_{kj}$ is strictly positive and defined such that
\begin{equation}
    \Lambda:= \Lambda^{kj} 
        = \{\lambda \mid \exp^k_{\lambda} \in L^1(\mathbb{R}^3) \}
        = \lbrace \lambda \in \mathbb{R}^5 \mid \lambda_2 <0 \rbrace
\end{equation}
is independent of $k$ and $j$.  
\end{ass} 
Roughly speaking, these assumptions are used to ensure integrability properties that are satisfied when the collision frequencies are independent of the velocity.  They are used in the technical details of the proofs below, but are in practice satisfied by many realistic frequency models.

\subsection{Target functions for intra-species collisions}
\label{sec:single_species_rigorou}

%\begin{equation}
%\mu^k(g) := \int a^k(v) g dv 
%\qquad
%\mu^k_*(g) := \mu^k(\nu_{kk}  g) = \int \nu_{kk} a^k(v) g dv 
%\end{equation}

%Associated to $g$ are moments $(n^k,u^k,T^k)$ and $(n^k_*,u^k_*,T^k_*)$ defined such that
%\begin{equation}
%\mu^k(g) = \begin{pmatrix} n^k \\ m_k n^k u^k \\ m_k n^k |u^k|^2 + 3 n^k T^k\end{pmatrix}
%\qquad 
%\mu^k_*(g) = \begin{pmatrix} n^k_* \\ m_k n^k_* u^k_* \\ m_k n^k_* |u^k_*|^2 + 3 n^k T^k_*\end{pmatrix}
%\end{equation}

%
%\subsubsection{Existence of solutions to the conservation constraints}
%\begin{lem}
%Assume that $\nu_{kk}>0$ and measurable on $\mathbb{R}^3$. Then, the set $\chi_k$ is non-empty.
%\end{lem}
%\begin{proof}
%Define a Maxwellian $M_k$ such that $\mu^k(M_k) = \mu^k_*(f_k)$
%\begin{equation}
%M_k(v) = n^k_* \left(\frac{m_k }{2 \pi T^k_*}\right)^{3/2} \exp\left(\frac{m_k|v-u^k_*|^2}{2 T^k_*} \right)
%\end{equation}
%Clearly 
%\begin{equation}
%\mu^k_*(M_k / \nu_{kk}) = \mu^k (M_k / \nu_{kk} ) = 
%\end{equation}
%Then $g=\frac{M_k}{\nu_{kk}}$ is a solution to \eqref{constraints2}.
%\end{proof}

%\subsubsection{Existence and Uniqueness of the target Maxwellians $M_{\lambda_{11}}$ and $ M_{\lambda_{22}}$}

%For example, if $\nu_{kk}$ is a polynomial in $v$, then we have $\Lambda^k=\lbrace \lambda \in \mathbb{R}^5 | \lambda_2 <0 \rbrace$, and thus $\Lambda^k\neq \emptyset$ and $\Lambda^k \bigcap \partial \Lambda^k= \lbrace \lambda \in \mathbb{R}^5 s.t. \lambda_2<0\rbrace \bigcap \lbrace \lambda \in \mathbb{R}^5 s.t. \lambda_2=0\rbrace= \emptyset$. 

We start the intra-species case; that is, for $ k \in \{1,2\}$, we show the existence of multiplier $\lambda^{kk}$ such that $M_{kk}$ satisfies \eqref{constraints_intra}.
The basic idea is to show that the dual function
\begin{align}
z(\lambda; \rho) = \mu^{kk}_0(\exp^k_{\lambda}) - \lambda \cdot \rho
\end{align}
is differentiable and attains its minimum on $\Lambda$ for any $\rho \in \mu^{kk}(D_{kk})$. Then the necessary condition for an extremum in $\Lambda$ yields 
\begin{align}
\label{eq:first-order-opt}
0= \nabla_{\lambda} z(\lambda^{kk}) = \int \nu_{kk}(v) a^k(v) \exp(\lambda^{kk} \cdot a^k(v)) dv - \rho,
\end{align}
which gives $\rho = \mu^{kk}(\exp^k_{\lambda^{kk}})$.  

\begin{lem}
\label{lem:Frechet}
The function $z$ is strictly convex and twice Fr\'echet differentiable on $\Lambda$. %Moreover, the Hessian is positive definite.
\end{lem}
\begin{proof}
It is sufficient to prove that  $\phi(\lambda) = \mu^{kk}_0(\exp^k_{\lambda})$ is strictly convex and twice Fr\'echet differentiable, with first derivative $D\phi(\lambda)= \mu^{kk}(\exp^k_{\lambda})$ and Hessian $H\phi(\lambda)= \int a^k(v) \otimes a^k(v) \exp_{\lambda}^k dv$. Convexity following immediately from convexity of the exponential function and linearity of the integral.  Specifically, given $\lambda^{(1)}$, $\lambda^{(2)}$ and two positive scalars $\theta_1, \theta_2$ such that $\theta_1 + \theta_2 = 1$, it follows that $\exp^k_{\theta_1 \lambda^{(1)} + \theta_2 \lambda^{(2)}} \leq \theta_1 \exp^k_{\lambda^{(1)}} + \theta_2 \exp^k_{\lambda^{(2)}}$.  Hence
\begin{align}
\begin{split}
    \phi(\theta_1 \lambda^{(1)} + \theta_2 \lambda^{(2)} ) 
    = \mu^{kk}_0(\exp^k_{\theta_1 \lambda^{(1)} + \theta_2 \lambda^{(2)}}) 
   & \leq \mu^{kk}_0(\theta_1 \exp^k_{\lambda^{(1)}} + \theta_2 \exp^k_{\lambda^{(2)}}) 
  \\ &  = \theta_1\phi( \lambda^{(1)}) + \theta_1 \phi(\lambda^{(2)} ).
\end{split}
\end{align}
For any nonzero $\delta \in \mathbb{R}^5$
\begin{align}
  \frac{  \phi (\lambda + \delta)- \phi(\lambda) - D \phi(\lambda) \cdot \delta}{|\delta|} = \int f_{\delta}(v) dv,
\end{align}
where
\begin{equation}
    f_{\delta}(v) = \nu_{kk}(v) \exp^k_\lambda (v) \left( \frac{\exp^k_\delta(v) - 1 - a^k(v) \cdot \delta}{|\delta|} \right).
\end{equation}
A Taylor series expansion shows that 
\begin{align}
\begin{split}
    \left| \frac{\exp^k_\delta(v) - 1 - \delta \cdot a^k(v)}{|\delta|}\right| = \left| \sum_{n=2}^{\infty} \frac{(\delta \cdot a^k(v))^n}{n!} \frac{1}{|\delta|} \right|& \leq |a^k(v)| \sum_{n=1}^{\infty} \frac{|\delta \cdot a^k(v)|^n}{n!} \\ &\leq |a^k(v)| \exp(|\delta \cdot a^k(v)|)
\end{split}
\end{align}
Therefore $f_{\delta}(v) \leq \exp^k_{\lambda/2}(v)  g_{\delta}(v)$, where
\begin{equation}
\begin{split}
    g_{\delta}(v)
       & := \nu_{kk}(v) |a^k(v)| \exp^k_{\lambda/2}(v) \exp(|\delta \cdot a^k(v)|) 
      \\&  \leq \nu_{kk}(v) |a^k(v)| \left( \exp^k_{\lambda/2+\delta}(v) + \exp^k_{\lambda/2-\delta}(v)  \right).
\end{split}
    %\nu_{kk}(v) \exp(\frac{\lambda}{2} \cdot a^k(v)) \exp (\frac{\lambda}{2} \cdot a^k(v)) |a^k(v)| \exp(\delta \cdot a^k(v)) =: \nu_{kk}(v) \exp(\frac{\lambda}{2} \cdot a^k(v)) g_{\delta}(v)
\end{equation}
Because $\Lambda$ is open, for $|\delta|$ sufficiently small, $\exp_{\lambda/2+\delta}(v)$ and $\exp_{\lambda/2-\delta}(v)$ are elements of $D^{kk}$, in which case $g_{\delta}$ is integrable.  Moreover, $\exp_{\lambda/2}$ is bounded. Hence $f_{\delta}$ is bounded above by an integrable function and the dominated convergence theorem gives
\begin{align}
    \lim_{\delta \rightarrow 0} \int f_{\delta}(v) dv = \int \lim_{\delta \rightarrow 0} f_{\delta}(v) dv =0.
    \end{align} 
    % Similar, for higher derivatives, lemma \ref{lem_est} yields integrable majorants for derivatives of $\exp^k_{\lambda}$. Thus $\lambda \mapsto \mu_0(\exp^k_{\lambda})\in C^{\infty}(\Lambda_k)$. \\
    
The existence of the Hessian can be proven in an analogous way.
\end{proof}

\begin{lem}\label{lem:min_line}
For fixed ${\lambda} \in \Lambda$, $\xi \in S^5$, and $\rho \in \mu^{kk}(D_{kk})$, the function
\begin{align}
    z_\xi (s) = z({\lambda} +s \xi;\rho)
    \label{function_s}
\end{align}
attains its unique minimum in the open interval
\begin{align}
    I(\xi, {\lambda}):= (- s_{\rm{b}}(- \xi, {\lambda}), s_{\rm{b}}(\xi, {\lambda}))
\end{align}
where
\begin{align*}
    s_{\rm{b}}(\xi, {\lambda}):= \sup \lbrace s: {\lambda} + s \xi \in \Lambda \rbrace
\end{align*}
takes the value $+ \infty$ if the boundary $\partial \Lambda$ is not met in the direction $\xi.$
\end{lem}

\begin{proof}
The fact that $z$ is strictly convex and differentiable with respect to $\lambda$ implies that $z_\xi$ is strictly convex and differentiable with respect to $s$.  Hence it attains a unique minimum on the closure of $I(\xi, {\lambda})$.  

We now show that $z_{\xi}$ cannot attain its minimum on the boundary of $I(\xi, {\lambda})$.  Suppose first that $s_{\rm{b}}(\xi,{\lambda}) < \infty$. 
According to Assumption \ref{ass:nu}, ${\lambda} + s_{\rm{b}}(\xi, {\lambda}) \xi \not \in \Lambda$.  Hence by Fatou's Lemma,
\begin{equation}
 \lim_{s \to s_{\rm{b}}(\xi, {\lambda})} \int \nu_{kk} \exp^k_{{\lambda} + s \xi} dv  
    \geq \int \nu_{kk}  \exp^k_{{\lambda} + s_{\rm{b}}(\xi, {\lambda}) \xi} dv = \infty
\end{equation}
which implies that  $ \lim_{s \to s_{\rm{b}}(\xi, {\lambda})}z_\xi(s) = + \infty$.

Suppose now that $s_{\rm{b}}(\xi, {\lambda})= \infty$.  There are two cases:
\begin{description}
    \item{\bf Case 1:} $\xi \cdot a^{k}(v) \leq 0$ for a.e. $v \in \mathbb{R}^3$.  Since $\rho \in \mu^{kk}(D_{kk})$, there exists $g \in D_{kk}$ such that $\rho = \mu^{kk}(g)$.  By definition,  $g$ is not identically zero and by Assumption \ref{ass:nu} $\nu_{kk}>0$. Thus the set
\begin{align}
    \Omega := \lbrace v \in \mathbb{R}^3 \mid \xi \cdot a^{k}(v) <0\rbrace \cap \lbrace v \in \mathbb{R}^3 \mid \nu_{kk}(v) g(v) >0 \rbrace
\end{align}
has positive measure.  Hence
    \begin{align}
    \xi \cdot \rho = \xi \cdot \mu^{kk}(g) = \int \nu_{kk}(v) \xi \cdot a^{k}(v) g(v) dv < 0
    \end{align}
    so that
    \begin{equation}
        \lim_{s \to \infty} z_{\xi}(s) 
            = \lim_{s \to \infty} \int \exp^{k}_{\lambda + s \xi} dv - (\lambda + s \xi) \cdot \rho 
            \geq \lim_{s \to \infty} - (\lambda + s \xi) \cdot \rho = \infty.
    \end{equation}
    \item{\bf Case 2:}  $\lbrace v \in \mathbb{R}^3: \xi \cdot a^k(v)>0 \rbrace$ has positive measure.  
    
    Then there exists an $\varepsilon>0$ such that $B=\lbrace v \in \mathbb{R}^3: \xi \cdot a^k(v) \geq \varepsilon\rbrace$ has positive measure. Hence 
\begin{align}
    \lim_{s \to \infty} z_{\xi}(s) \geq \lim_{s \to \infty} \left( \left(\int_B \nu_{kk}(v) \exp^k_{{\lambda}} dv \right) \exp (s \varepsilon) - ({\lambda}+ s \xi) \rho \right) = \infty
\end{align}
 due to exponential growth in $s$.
\end{description}

% We conclude that on each endpoint of the interval $I(\xi, {\lambda})$, $z(s)$ either diverges or the slope becomes $\infty$. In any case, there exist points $s_1,s_2 \in I(\xi, {\lambda})$ such that $z'(s_1) z'(s) <0$. Applying the mean value theorem, we find $s^* \in I(\xi, {\lambda})$ such that $z'(s^*)=0$ and due to strict convexity of $z$, $s^*$ is the unique minimizer.

\end{proof}
\begin{thm}\label{thm:min}
% Assume $\emptyset \neq \Lambda \subset \mathbb{R}^{N+1}$ is open and convex. Further, let $z:\Lambda \mapsto \mathbb{R}$ be strictly convex and twice differentiable. If for ${\lambda} \in \Lambda$
% \begin{align}
% \begin{split}
% s \mapsto z({\lambda} + s \xi), \quad s \in I(\xi, {\lambda}):= ( - s_{\rm{b}}(- \xi, {\lambda}), s_{\rm{b}}(\xi, {\lambda} ), \\ s_{\rm{b}}(\xi, {\lambda}):= \sup \lbrace s: {\lambda} + s \xi \in \Lambda \rbrace
% \end{split}
% \end{align} 
% attains its minimum for all directions $\xi \in S^N$, then 
For any $\rho \in \mu^{kk}(D_{kk})$, the function $z(\cdot;\rho)$ has a unique minimizer $\lambda^* \in \Lambda$.
\end{thm}
\begin{proof}
Let $\lbrace \lambda^{(\ell)} \rbrace_{\ell=0}^{\infty}$ be an infimizing sequence such that $z(\lambda^{(\ell)}) \rightarrow z_*$, where 
\begin{align*}
    z_* = \inf_{\lambda \in \Lambda} z(\lambda).
\end{align*}
Let $d^{(\ell)}=\lambda^{(\ell)}-\lambda^{(0)}; ~ \ell \geq 1$ and set $\xi^{(\ell)} = {d^{(\ell)}}/{||d^{(\ell)}||}$. Then $\xi^{(\ell)} \rightarrow \xi^* \in S^{4}$ possibly via a subsequence, because $S^{4}$ is compact. For any $\xi \in S^{4},$ let $s_*(\xi)= \arg \min_{s \in \mathbb{R}} z(\lambda^{(0)} + s \xi; \rho)$ which, according to Lemma \ref{lem:min_line}, is well-defined. Because $z$ is strictly convex and twice differentiable, 
\begin{align*}
  &(i) \quad g(\xi, s):= \partial_s z(\lambda^{(0)} + s \xi; \rho)=0 \quad \text{if and only if} \quad s=s_*(\xi) \\ 
    &(ii) \quad \partial_s g(\xi, s) >0 
\end{align*}
Thus the implicit function theorem implies that $s_*$ is a $C^1$ function in a neighbourhood $N(\xi_*) \subset \Lambda$ that satisfies 
\begin{align}
    g(\xi, s_*(\xi)) =0.
\end{align}
Let $\ell_*$ be large enough that
   $ \xi^{(\ell)} \in N(\xi_*) $ for all $\ell \geq \ell_*$. Then
   \begin{align}
       z(\lambda^{(\ell)}; \rho) 
        = z(\lambda^{(0)} + d^{(\ell)}; \rho) 
        = z(\lambda^{(0)} + ||d^{(\ell)}|| \xi^{(\ell)}; \rho) 
        \geq z(\lambda^{(0)} + s_*(\xi^{(\ell)}) \xi^{(\ell)}; \rho).
       \label{min_line}
   \end{align}
   Because $s_*$ is continuous on $N(\xi_*)$ the sequence $s_*(\xi^{(\ell)}) \rightarrow s_*(\xi^{*})$ with $|s_*(\xi^{*})|< \infty$.  Moreover, since $z$ is continuous
   \begin{align}
       z_* = \lim_{\ell \rightarrow \infty}
z(\lambda^{(\ell)}) \geq \lim_{\ell \rightarrow \infty} z(\lambda^{(0)} + s_*(\xi^{(\ell)}) \xi^{(\ell)}) = z(\lambda^{(0)} + s_*(\xi^{*})  \xi^*) \geq z_*,   \end{align}
where first inequality follows from \eqref{min_line}. Hence the infimum is attained at $\lambda_* = \lambda^{(0)} + s_*(\xi^{*})  \xi_* \in \Lambda.$
\end{proof}

\begin{corollary}
Given any $f_k \in D_{kk}$, there exists a unique multiplier $\lambda^{kk}$ such that $M^{kk}$ given by \eqref{form1} solves \eqref{min}.
\end{corollary}
\begin{proof}
Let $\rho_k = \mu^{kk}(f_k)$.  According to Theorem \ref{thm:min}, $z(\cdot,\rho_k)$ has a unique minimizer in $\Lambda$, which we denote by $\lambda^{kk}$.  By Lemma \ref{lem:Frechet}, 
$z(\cdot,\rho_k)$ is also differentiable, so the first-order optimality condition \eqref{eq:first-order-opt} implies that $\rho_k = \mu^{kk}(\exp_{\lambda^{kk}})$.  The result then follows from Theorem \ref{thm:unique_min}.
\end{proof}

\subsection{Target functions for inter-species collisions}
In this section we show the existence of the multipliers $\lambda^{12} = (\lambda_0^{12}, \lambda_1^{12}, \lambda_2^{12}) \in \mathbb{R} \times \mathbb{R}^3 \times \mathbb{R}$  and $\lambda^{21} =(\lambda_0^{21}, \lambda_1^{21}, \lambda_2^{21}) \in \mathbb{R} \times \mathbb{R}^3 \times \mathbb{R}$ such that $\lambda_1^{12}=\lambda_1^{21}$, $\lambda_2^{12}=\lambda_2^{21}$, and $M_{12}$ and $M_{21}$ satisfy \eqref{constraints_inter}. Denote
\begin{align}
{\lambda} =(\lambda^{1}_0, \lambda^{2}_0, \lambda_1, \lambda_2) 
\quad {\lambda}^{1} = (\lambda^{1}_0, \lambda_1, \lambda_2) 
\quad \lambda^{2} = (\lambda^{2}_0, \lambda_1, \lambda_2) 
\end{align}
and use this notation for other vectors when appropriate. 
%
% If we consider for example a $\xi\in \mathbb{R}^3$, then $\xi^{12}$ is a vector in $\mathbb{R}^5$ which is obtained by taking $\xi$ without the second component $\xi_0^{21}$, whereas $\xi^{21} \in \mathbb{R}^5$ is obtained by taking $\xi$ without the first component $\xi_0^{12}$.
%
 % \\
%\tilde{\lambda}:=(\tilde{\lambda}_0^{12}, \tilde{\lambda}_0^{21}, \tilde{\lambda_1}, \tilde{\lambda}_4):= - \lambda,~ \tilde{\lambda}_{1,3-6} = - \lambda_{1,3-6},~ \tilde{\lambda}_{2-6}=- \lambda_{2-6}
% We further let
% \begin{align}
% \exp_{\lambda^{12}}(v):= \exp (\lambda^{12} \cdot a^1(v))
% \quad
% \text{and}
% \quad
% \exp_{\lambda^{21}}(v):= \exp (\lambda^{21} \cdot a^2(v)) 
% \end{align}
Given $g_1,g_2 \in D$, let
\begin{align}
\bar{\mu}(g_1,g_2) = 
\begin{pmatrix}
\mu_0^{12}(g_1) \\
\mu_0^{21}(g_2) \\
\mu_1^{12}(g_1) + \mu_1^{21}(g_2) \\
\mu_2^{12}(g_1) + \mu_2^{21}(g_2)
\end{pmatrix}.
% $ ( \mu_0^{12}(g_1) , 0 , \mu_1^{12}(g_1) , \mu_2^{12}(g_1))^T + ( 0, \mu_0^{21}(g_2), \mu_1^{21}(g_2), \mu_2^{21}(g_2))^T %\mu_0^{12}(g_1) +  \mu_0^{12}(g_2)%:= \int \nu_{12} g_1 dv +  \int \nu_{21} g_2 dv $
\end{align}
%
%We make the following assumption on the collision frequencies.
%\begin{ass}
%\label{ass12}
%We assume that $\nu_{21}$ has the following property: If $g$ is a function such that $\nu_{12}(v)g(v)$ is integrable, then $\nu_{21}(v)g(v)$ should be also integrable.
%\textcolor{cyan}{It might be also possible without this assumption, but it becomes much more technical so I first tried it with this assumption.}
%\end{ass}
%Next, we define the set 
%$$D=\lbrace f \geq 0: \nu_{12} (1+ |v|^2) f \in L^1(\Omega), f \neq 0 \rbrace,$$
%the  admissible set of parameters, 
%\begin{align}
%\Lambda := \lbrace \lambda \in \mathbb{R}^5: \exp^k_{\lambda} \in D \rbrace
%\end{align}
%the corresponding exponential functions,
%\begin{align}
%E:= \lbrace \exp^k_{\lambda} : \lambda \in \Lambda \rbrace,
%\end{align}
%and introduce the notation 
For any $\bar \rho \in \bar{\mu}(D_{12} \times D_{21})$, introduce the dual function
\begin{align}
\label{eq:zmix}
\bar{z}(\lambda; \bar \rho) =  \mu^{12}_0(\exp^1_{\lambda^1}) +  \mu^{21}_0(\exp^2_{\lambda^2})  - \lambda \cdot \bar \rho. %\\ &\rho_{\rm{\rm{mix}}}:= \mu^{\rm{\rm{mix}}}(g_1,g_2)
\end{align}
Similar to the intra-species case, our goal is to show that for any such $\bar \rho$, $z(\lambda; \bar \rho)$ attains its minimum on 
\begin{equation}
    \bar{\Lambda}= \lbrace \lambda \in \mathbb{R}^6: \lambda^{1}, \lambda^{2} \in \Lambda \rbrace.
\end{equation}
Then the necessary first-order condition for a minimum at $\lambda$ 
\begin{align}
0= \nabla_{\lambda} z(\lambda;\bar \rho) 
= \bar{\mu}({\exp^1_{\lambda^{1}}(v), \exp^2_{\lambda^{2}}(v)})
- \bar\rho,
\label{eq:first-order-opt_inter}
\end{align}
which recovers the required constraints in \eqref{constraints_inter}, if we set $\lambda^{12}=\lambda^{1}$ and $\lambda^{21}=\lambda^{2}$.

%We also make the following assumption.
%\begin{ass}
%Assume that $\nu_{12}$ is such that $\Lambda \neq \emptyset$ and $\Lambda \bigcap \partial \Lambda = \emptyset$. 
%\label{ass_coll2}
%\end{ass}
%We make the following general assumptions
%With these definitions and  assumptions it is easy to see that Lemma \ref{lem_est} holds for the collision frequencies $\nu_{12}$ and $\nu_{21}$ and the exponential functions $\exp_{\lambda_{1,3,4,5,6}}$ and $\exp_{\lambda_{2,3,4,5,6}}$, if we make the same assumptions on $\nu_{12}$ and $\nu_{21}$ as are stated in Lemma \ref{lem_est}. Thus, we will just state the lemma without giving the proof.
\begin{lem}
\label{lem3.2.2}
The function $\bar{z}$ defined in \eqref{eq:zmix} is strictly convex and twice Fr\'echet differentiable on $\bar \Lambda$.
\label{lem_Frechet2}
\end{lem}
\begin{proof}
Differentiability of the $\bar z$ can be deduced as in the intra-species case by simply following the arguments of Lemma \ref{lem:Frechet}.  We skip these details.   Convexity also follows in a similar way.  Let $\bar{\phi}(\lambda) = \mu^{12}_0(\exp^1_{\lambda^1}) +  \mu^{21}_0(\exp^2_{\lambda^2})$, then convexity of the exponential function implies that for any $\theta \in (0,1)$, $\lambda \in \bar{\Lambda}$, and $\beta \in \bar {\Lambda}$,
\begin{align}
\begin{split}
    \bar{\phi}(\theta \lambda) + \bar{\phi}((1-\theta) \beta ) 
    &= \mu^{12}_0(\exp^1_{\theta \lambda^1 + (1-\theta) \beta ^1}) + \mu^{21}_0(\exp^2_{\theta \lambda^2 + (1-\theta) \beta ^2}) \\
    &\leq \mu^{12}_0(\theta \exp^1_{\lambda^1} + (1-\theta) \exp^1_{\beta^1}) + \mu^{21}_0(\theta \exp^2_{\lambda^2} + (1-\theta) \exp^2_{\beta^2})  \\
    &= \theta \bar{\phi}(\lambda) + (1-\theta) \bar{\phi}(\beta ) 
\end{split}
\end{align}
Thus $\bar{\phi}$ is strictly convex, as is $\bar{z}$, since the two functions differ only by a linear term.
\end{proof}
\begin{lem}
For ${\lambda} \in \bar \Lambda , \xi \in S^5$, and $\bar \rho \in \bar\mu(D_{12} \times D_{21})$, the function 
\begin{align}
  \bar z_{\xi} \colon   s \mapsto \bar z({\lambda}+ s \xi; \bar \rho)
\end{align}
attains its unique minimum in the open interval 
\begin{align}
    \bar I(\xi, {\lambda}):= (- \bar s_{\rm{b}}(- \xi, {\lambda}), \bar s_{\rm{b}}(\xi, {\lambda})), 
\end{align}
where 
\begin{align}
   \bar s_{\rm{b}}(\xi, \lambda)= \sup \lbrace s \colon \lambda^1 + s \xi^1, \lambda^2 + s \xi^2 \in \Lambda \rbrace.
\end{align}
\end{lem}
\begin{proof}
We follow the arguments of the proof of Lemma \ref{lem:min_line}.  The fact that $\bar z$ is strictly convex and differentiable with respect to $\lambda$ implies that $\bar z_\xi$ is strictly convex and differentiable with respect to $s$.  Hence $\bar z_\xi$ attains a unique minimum on the closure of $\bar I(\xi, {\lambda})$. 
We therefore need only show that $\bar z_{\xi}$ cannot attain its minimum on the boundary of $\bar{I}(\xi, {\lambda})$. 

Suppose first that $\bar{s}_{\rm{b}}(\xi, \lambda) < \infty$. By Fatou's Lemma,
%\begin{align}
%z'(s) &= \left( \begin{pmatrix}
%\mu_0^{12}(\exp_{\lambda^{12}+s \xi}) \\ 0 \\  \mu_1^{12}( \exp_{\lambda^{12}+ s \xi}) \\  \mu_2^{12}( \exp_{\lambda^{12}+ s \xi})
%\end{pmatrix} - \begin{pmatrix} \mu_0^{12}(f_1) \\ 0 \\ \mu_1^{12}(f_1) \\ \mu_2^{12} (f_1) \end{pmatrix} \right) \cdot \xi -  \left( \begin{pmatrix}
%0 \\ \mu_0^{21}(\exp_{\lambda_{21}+s \xi}) \\  \mu_1^{21}( \exp_{\lambda_{21}+ s \xi}) \\  \mu_2^{21}( \exp_{\lambda_{21}+ s \xi})
%\end{pmatrix} - \begin{pmatrix} 0 \\ \mu_0^{21}(f_2) \\ \mu_1^{21}(f_2) \\ \mu_2^{21}(f_2) \end{pmatrix} \right) \cdot \xi \rightarrow \infty
%\end{align}
%for $s \rightarrow s_{\rm{b}}^{\rm{\rm{mix}}}(\xi, {\lambda}).$
%This is because if $s \rightarrow s_{\rm{b}}(\xi,\lambda)$, then

\begin{multline}
 \lim_{s \to \bar s_{\rm{b}}(\xi, {\lambda})} 
 \left\{\int \nu_{12} \exp^1_{{\lambda^1} + s \xi^1} dv  + \int \nu_{21}\exp^2_{{\lambda^2} + s \xi^2} dv  \right\} \\
    \geq \left\{\int \nu_{12} \exp^1_{{\lambda^1} + \bar s_{\rm{b}}(\xi, {\lambda}) \xi^1} dv  + \int \nu_{21}\exp^2_{{\lambda^2} + s_{\rm{b}}(\xi, {\lambda}) \xi^2} dv  \right\} dv 
\end{multline}
Assumption \ref{ass:nu} implies that ${\lambda^1} + \bar s_{\rm{b}}(\xi, {\lambda}) \xi^1 \not \in \Lambda$ or ${\lambda^1} + \bar s_{\rm{b}}(\xi, {\lambda}) \xi^1 \not \in \Lambda$.  Hence at least one of the integrals on the right-hand side above is $\infty$, which implies
\begin{align}
    \lim_{s \to \bar s_{\rm{b}}(\xi, \lambda)}  z_{\xi}(s)
        =  \mu^{12}_0(\exp^1_{\lambda^1}) +  \mu^{21}_0(\exp^2_{\lambda^2})  - \lambda \cdot \bar \rho
        = \infty.
\end{align}
%Therefore  the first or the second term goes to infinity which can be proven as in the proof of lemma \ref{lem:min_line} exchanging the intra species collision frequencies by the inter species collision frequencies. 

Now suppose instead that $\bar s_{\rm{b}}(\xi, \lambda)= \infty$. There are two cases:

\noindent {\bf Case 1:} %$s_{\rm{b}}(\xi, {\lambda})=\infty$. %We introduce 
%\begin{align}
%L^+:= \lbrace \xi \in S^6: \xi \cdot \rho \geq 0 \rbrace \quad \text{and} \quad L^-:= \lbrace \xi \in S^6: \xi \cdot \rho <0\rbrace.
%\end{align}
%For $\xi \cdot \rho <0$, we have $-s \xi \cdot \rho \rightarrow \infty$ for $s \rightarrow s_{\rm{b}}(\xi, \lambda)=\infty$.
%Since $\mu_0$ is always positive, we see that in this case $z(s) \rightarrow \infty.$
$\xi^{1} \cdot a^1(v) \leq 0$ and $\xi^{2} \cdot a^2(v) \leq 0$ for a.e $v \in \mathbb{R}^3$. 

Since $\bar \rho \in \bar{\mu} (D_{12} \times D_{21})$, there exist $g_1, g_2 \in D_{12} \times D_{21}$ such that $\bar{\rho} = \bar{\mu}(g_1,g_2)$;  that is 

\begin{align}
\bar{\rho} = \bar{\mu}(g_1,g_2) = 
\begin{pmatrix}
\mu_0^{12}(g_1) \\
\mu_0^{21}(g_2) \\
\mu_1^{12}(g_1) + \mu_1^{21}(g_2) \\
\mu_2^{12}(g_1) + \mu_2^{21}(g_2)
\end{pmatrix}.
\end{align}
By definition,  $g_1$ and $g_2$ are not identically zero, and by Assumption \ref{ass:nu}, $\nu_{kj}>0$. Thus the sets
\begin{align}
    \Omega_1 &:= \lbrace v \in \mathbb{R}^3 \mid \xi^{1}  \cdot a^1(v)<0\rbrace \cap \lbrace v \in \mathbb{R}^3  \mid \nu_{12}(v) g_1(v) >0 
    \rbrace  \quad \text{and}\\
    \Omega_2 &:= \lbrace v \in \mathbb{R}^3 \mid \xi^{2} \cdot a^2(v) <0 \rbrace \cap \lbrace v \in \mathbb{R}^3  \mid \nu_{21}(v) g_2(v) >0 
    \rbrace 
\end{align}
both have positive measure. Hence
\begin{align}
\xi \cdot \bar{\rho} 
    &= \xi^{1} \cdot  \mu^{12}(g_1) +  \xi^{2} \cdot \mu^{21}(g_2)\\ 
    &= \int \nu_{12} \xi^{1} \cdot a^1(v) g_1(v)dv +  \int \nu_{21} \xi^{2} \cdot a^2(v) g_2(v) dv <0,
\end{align}
so that
\begin{align}
    \lim_{s \rightarrow \infty} \bar z_{\xi}(s) 
    &=  \lim_{s \rightarrow \infty}  \left\{ \mu^{12}_0(\exp^1_{\lambda^1+s \xi^1}) +  \mu^{21}_0(\exp^2_{\lambda^2 + s \xi^2})  - (\lambda + s \xi)\cdot \bar \rho \right \} \\
    &> \lim_{s \rightarrow \infty} \left\{- (\lambda + s \xi) \cdot \bar{\rho} \right\} = \infty.
\end{align}

\noindent {\bf Case 2:} 
The set $\lbrace v \in \mathbb{R}^3 \mid \xi^{1} \cdot a^1(v) >0\rbrace $ or $\lbrace v \in \Omega \mid \xi^{2} \cdot a^2(v)>0 \rbrace $ has positive measure. 

Without loss of generality, assume that $\lbrace v \in \mathbb{R}^3 \mid \xi^{1} \cdot a^1(v) >0\rbrace $ has positive measure. Then, there exists some $\varepsilon >0$ such that  $B=\lbrace v \in \mathbb{R}^3 \mid \xi^{12}\cdot a^1(v)> \varepsilon \rbrace$ also has positive measure. Hence
\begin{align}
\lim_{s \to \infty} \bar{z}_{\xi}(s) \geq \lim_{s \to \infty} \left( \left( \int_B \nu_{12} \exp^1_{\lambda^{1}} dx \right) \exp (s \varepsilon) - ({\lambda} + s \xi) \cdot \rho_{\rm{\rm{mix}}} \right) 
= \infty.
\end{align} 
   due to exponential growth in $s$.

%We conclude that on each endpoint of the interval $I(\xi, \lambda)$, $z_{\rm{\rm{mix}}}(s)$ either diverges or that the slope $z'_{\rm{\rm{mix}}}(s)$ becomes $+ \infty$. In either case, there exist points $s_1,s_2 \in I(\xi, \lambda)$ such that $z'_{\rm{\rm{mix}}}(s_1) z'_{\rm{\rm{mix}}}(s_2)<0$. Applying the mean value theorem, we can find $s^*\in I(\xi, \lambda)$ such that $z'_{\rm{\rm{mix}}}(s^*)=0$ and, due to strict convexity of $z_{\rm{\rm{mix}}}$, $s^*$ is the unique minimizer. %The rest of the proof is very similar to the proof of Theorem 2 in \cite{Junk}, so we will skip it here.
\end{proof}
\begin{thm}
\label{thm:unique_min_inter}
For any $\bar \rho \in \bar \mu(D_{12} \times D_{21})$, the function $\bar{z}(\cdot, \bar \rho)$ has a unique minimizer $\lambda^* \in \bar \Lambda$.
\end{thm}
The proof of this theorem is analogous to the proof of Theorem \ref{thm:min} in the intra-species case.  
\begin{corollary}
Given any $f_1 \in D_{12}$ and $f_2 \in D_{21}$, there exist multipliers $\lambda^{12}$ and $\lambda^{21}$ such that $\lambda^{21}_1 = \lambda^{12}_1$,  $\lambda^{21}_2 = \lambda^{12}_2$, and the corresponding functions $M^{12}$ and $M^{21}$ given in \eqref{form1} solve \eqref{constraints_inter}.
\end{corollary}
\begin{proof}
Let $\bar \rho= \bar \mu(f_1,f_2)$. According to Theorem \ref{thm:unique_min_inter}, $\bar z(\cdot, \bar \rho)$ has a unique minimizer, which we denote by $\lambda^*=((\lambda^*)_0^{1},(\lambda^*)_0^{2}, (\lambda^*)_1, (\lambda^*)_2).$ By Lemma \ref{lem:Frechet}, $\bar z(\cdot, \bar \rho)$ is also differentiable, so the first-order optimality condition \eqref{eq:first-order-opt_inter} implies that $\bar \rho=  \bar{\mu}(\exp^1_{(\lambda^*)^{1}}, \exp^2_{(\lambda^*)^{2}})$. The result then follows from Theorem \ref{thm2.2}. Finally, we set 
\begin{equation}
    \lambda^{12} = ((\lambda^*)_0^{1}, (\lambda^*)_1, (\lambda^*)_2) 
    \quad \text{and} \quad 
    \lambda^{21} = ((\lambda^*)_0^{2}, (\lambda^*)_1, (\lambda^*)_2) 
\end{equation}
and define $M^{12}$ and $M^{21}$ according to \eqref{form1}.
\end{proof}
\section{Consistency of the model}
\label{sec4}
The conditions \eqref{constraints_intra} and \eqref{constraints_inter} lead to standard conservation laws and an entropy dissipation statement.  We recall a few definitions:
\begin{defn}
The mass density, momentum, and energy of an integrable distribution $g = g(v)$ of particles with mass $m$ are given by the moments
\begin{equation}
    \rho_g = \int m g(v) dv, \quad q_g = \int m v g(v) dv, \quad \text{and} \quad E_g = \frac12 \int m |v|^2 g(v) dv,
\end{equation}
respectively.  The associated mean velocity and temperature are given by
\begin{equation}
    u_g = \frac{q_g}{\rho_g} = \frac{\int v g(v) dv}{\int g(v) dv} 
        \quad \text{and} \quad 
    T_g = \frac{2}{3}\frac{E_g}{\rho_g/m} - \frac13 \frac{|q_g|^2}{\rho_g} = \frac{1}{3}\frac{\int m |v-u_g|^2 g(v) dv}{\int g(v) dv}.
\end{equation}
\end{defn}
\subsection{Conservation properties}
An immediate consequence of \eqref{constraints_intra} and \eqref{constraints_inter} is the following.
\begin{thm}[Conservation of the number of each species, total momentum and total energy]  The space-homogeneous form of \eqref{BGK2} satisfies
\begin{gather}
    \partial_t \rho_{f_1} = \partial_t \rho_{f_2}   =0, \quad 
    \partial_t \left( q_{f_1} + q_{f_2}  \right) = 0 , \quad
    \partial_t \left( E_{f_1} + E_{f_2}  \right) = 0 
\end{gather}

% Using the exponential form \eqref{form1} for the target Maxwellians, we have:
% \begin{itemize}
% \item conservation of the number of particles, i.e.,
% $$
% \int  Q_{11}(f_1,f_1) dv  =  \int  Q_{12}(f_1,f_2) dv  = \int Q_{22}(f_2, f_2) dv = \int Q_{21}(f_2,f_1) dv = 0,
% $$
% \item conservation of total momentum,
% $$
% \int m_1 v [Q_{11}(f_1,f_1)+Q_{12}(f_1,f_2)] dv +
% \int m_2 v [Q_{22}(f_2,f_2)+Q_{21}(f_2,f_1)] dv= 0,
% $$
% \item and conservation of total energy
% \begin{align}
% \int \frac{m_1}{2} |v|^2  (Q_{11}(f_1,f_1)+Q_{12}(f_1,f_2)) dv  +
% \int \frac{m_2}{2} |v|^2  (Q_{22}(f_2,f_2)+Q_{21}(f_2,f_1)) dv = 0.
% \end{align}
% \end{itemize}
 \end{thm}
\subsection{Entropy dissipation and the structure of equilibria}
Define the total entropy density
\begin{equation}
    H(g_1,g_2) = \int h(g_1)dv + \int h(g_2) dv
\end{equation}
and the dissipation density
\begin{align}
S(g_1,g_2) &= S_{11}(g_1) + S_{12}(g_1,g_2) + S_{21}(g_1,g_2) + S_{22}(g_2) \\
&=\int \nu_{11} \ln g_1 ( M_{11} - g_1) dv
+ \int \nu_{12} \ln g_1 ( M_{12}- g_1) dv \\
&\quad + \int \nu_{21} \ln g_2 ( M_{21}- f_2) dv 
+ \int \nu_{22} \ln g_2 ( M_{22}- g_2) dv
\end{align}
\begin{thm}
Assume $g_1,g_2>0$. Then $S(g_1,g_2) \geq 0$ with equality if and only if $g_1$ and $g_2$ are two Maxwellian distributions with equal mean velocity and temperature.
\end{thm}
\begin{proof}
In \cite{Struchtrup}, it is shown that $S_{kk}(g) \geq 0$ with equality if and only if $g$ is a Maxwellian. Thus it remains to show a similar result for the combined quantity $S_{12}(g_1,g_2) + S_{21}(g_1,g_2)$. 
% \begin{align}
% \int \nu_{12} \ln g_1 ( M_{12}- g_1) dv + \int \nu_{21}\ln g_2 (M_{21}- g_2) dv \leq 0.
% \end{align}
We begin with the following claim:
\begin{align}
\label{eq:I_def}
I(g_1,g_2):= \int \nu_{12} \ln M_{12} ( M_{12} - g_1) dv +  \int \nu_{21} \ln M_{21} ( M_{21} - g_2) dv = 0.
\end{align}
Indeed an explicit calculation gives 
\begin{align}
\ln M_{12}  =   m_1 \lambda_0^{12}  + m_1 \lambda_1\cdot v  + m_1 \lambda_2 |v|^2  \quad \text{and} \quad   
\ln M_{21}  =    m_2 \lambda_0^{21} + m_2 \lambda_1 \cdot v + m_2 \lambda_2 |v|^2,
\end{align}
which when substituted into \eqref{eq:I_def} gives
\begin{align}
I(g_1,g_2) 
    &= \int \nu_{12} ( m_1 \lambda_0^{12}+ m_1 \lambda_1\cdot v + m_1 \lambda_2 |v|^2) ( M_{12} - g_1) dv \\
    & \quad  +  \int \nu_{21} ( m_2 \lambda_0^{21}+ m_2 \lambda_1\cdot v + m_2 \lambda_2 |v|^2) ( M_{21} - g_2) dv = 0,
\end{align}
due to the constraints \eqref{constraints_inter}. From \eqref{eq:I_def}, it follows that
\begin{align}
\begin{split}
S_{12}(g_1,g_2) &+ S_{21}(g_1,g_2)
=  S_{12}(g_1,g_2) + S_{21}(g_1,g_2) - I(g_1,g_2) \\
&=  \int \nu_{12} \ln \left( \frac{g_1}{M_{12}} \right) ( M_{12} - g_1) dv 
    +  \int \nu_{21} \ln \left( \frac{g_2}{M_{21}} \right) ( M_{21} -g_2) dv \\
    &\leq 0.
\label{equality}
\end{split}
\end{align}
with equality if and only if $g_1=M_{12}$ and $g_2=M_{21}$.
Moreover, a direct calculation shows that the functions $M_{12}$ and $M_{21}$ have the same mean velocity and temperature:
% \begin{align}
% n_1:= \int g_1 dv = \frac{\sqrt{2 \pi}^3}{\sqrt{m_1 \lambda_2}} \exp \left(\frac{\lambda_1^2 m_1}{2 \lambda_2} - \lambda_0^{12}\right)
% \quad n_2:= \int g_2 dv = \frac{\sqrt{2 \pi}^3}{\sqrt{m_2 \lambda_2}} \exp \left(\frac{\lambda_1^2 m_2}{2 \lambda_2} - \lambda_0^{21} \right).
% \end{align}
\begin{align}
%\frac{\int m_1 v M_{12} dv }{\int M_{12} dv} = \frac{\int m_2 v M_{21} dv }{\int M_{21} dv} = \frac{\lambda_1}{\lambda_2}
u_{M_{12}} = u_{M_{21}} =- \frac{\lambda_1}{\lambda_2}
\quad \text{and} \quad
T_{M_{12}} = T_{M_{21}} =- \frac{1}{2 \lambda_2}
\label{equ_vel}
\end{align}

% \begin{align}
% |u_1|^2 + 3 \frac{T_1}{m_1}:= \frac{1}{n_1} \int |v|^2 g_1 dv = 3 \frac{1}{m_1 \lambda_2} + |\frac{\lambda_1}{\lambda_2}|^2, \\ |u_2|^2 + 3 \frac{T_2}{m_2}:= \frac{1}{n_2} \int |v|^2 g_2 dv = 3 \frac{1}{m_2 \lambda_2} + |\frac{\lambda_1}{\lambda_2}|^2.
% \end{align}
% Using \eqref{equ_vel}, we obtain that $T_1=T_2$.

% Therefore, we obtain two Maxwell distributions with equal mean velocities and temperatures.
\end{proof}

 \begin{corollary}[Entropy inequality for mixtures]
Assume that $f_1, f_2 >0$ are a solution to \eqref{BGK2} where the target Maxwellians have the shape  \eqref{form1}, then we have the following entropy inequality
\begin{align}
\partial_t \left(  H(f_1,f_2) \right) + \nabla_x \cdot \left(\int  v (h(f_1) + h(f_2)) dv \right) \leq 0
\end{align}
with equality if and only if $f_1$ and $f_2$ are two Maxwellian distributions with equal mean velocity and temperature.
\end{corollary}
\begin{proof}
A direct calculation with $\eqref{BGK2}$ gives
\begin{equation}
    \partial_t H(f_1,f_2) + \nabla_x \cdot \int (h(f_1)+ h(f_2)) v dv = S(f_1,f_2).
\end{equation}
The result then follows immediately from the previous theorem.
\end{proof}
\section{The \textit{N}-species case}
\label{sec5}
 The two-species case can be  extended to a system of $N$-species that undergo binary collisions. We consider the $N$-species kinetic equation,
\begin{align}
\partial_t f_i + v \cdot \nabla_x f_i = \sum_{j=1}^N \nu_{ij} (M_{ij}- f_i), \quad i=1,...,N.
\end{align}
%For future reference, we denote $Q(f_i,f_j)=\nu_{ij} (M_{ij}- f_i)$. 
The quantity $\nu_{ii} $ is the collision frequency of particles of species $i$ with itself whereas $\nu_{ij} $ is the collision frequency of particles of species $i$ with species $j$, with $i,j=1,...,N, ~ i \neq j$. %In the following, if we write $\nu_{ij}(v)$, we mean both the explicit dependence of $\nu_{ij}$ with respect to $v$ and the implicit dependence of $v$ coming from the dependence of $f_i$ and $f_j$.
We only have terms of this form and not terms containing indices of more than two species because we consider only binary interactions.

For fixed $i,j\in \{1,\dots,N\}$ the  target Maxwellians $M_{ii}$, $M_{jj}$, %  are defined as
%\begin{align}
%M_{ii} =  e^{\lambda_0^{ii} +\lambda_1^{ii} \cdot v +\lambda_2^{ii} |v|^2}, \quad i=1,...,N
%\label{oneN}
%\end{align}
%with coefficients $\lambda_0^{ii} \in \mathbb{R}, \lambda_1^{ii} \in \mathbb{R}^3, \lambda_2^{ii} \in \mathbb{R}^+$, $i=1,...,N.$
 $M_{ij}$ and $M_{ji}$ 
 are given by \eqref{form1}.
%\begin{align}
%M_{ij} =  e^{\lambda_0^{ij} + m_i \lambda_1^{ij} \cdot v + m_i \lambda_2^{ij} |v|^2}, \quad M_{ji} =  e^{\lambda_0^{ji} + m_j \lambda_1^{ji} \cdot v + m_j \lambda_2^{ji} |v|^2},
%\label{mixN}
%\end{align}
%with coefficients $\lambda_0^{ij} \in \mathbb{R}, \lambda_1^{ij}=\lambda_1^{ji} \in \mathbb{R}^3, \lambda_2^{ij}=\lambda_2^{ji} \in \mathbb{R}^+$, $i,j=1,...,N, ~ i \neq j$ possibly dependent on $x$ and $t$. 
The single species target Maxwellians $M_{ii}$ and $M_{jj}$  will be determined  such that they satisfy \eqref{constraints_intra}.
%\begin{align}
%\int \nu_{ii} a^i(v) (M_{ii} - f_i) dv = 0, \quad i=1,...,N. 
%\end{align}
%This guarantees conservation of the number of particles, momentum and energy in interactions of a particle with a particle of the same species.
The functions $M_{ij}$ and $M_{ji}$ will be determined such that we obtain conservation of mass of each species and conservation of total momentum and total energy in interactions between these two species, i.e.,
\begin{align}
\begin{split}
\int \nu_{ij} M_{ij} dv = \int \nu_{ij} f_i dv, \quad \int \nu_{ji} M_{ji} dv = \int \nu_{ji} f_j dv \\
  \int \nu_{ij} \begin{pmatrix}
m_i v \\ m_i |v|^2
\end{pmatrix} (M_{ij} - f_i) dv  = - 
\int \nu_{ji} \begin{pmatrix}
m_j v \\ m_j |v|^2
\end{pmatrix} (M_{ji} - f_j) dv.
\end{split}
\end{align}
as an obvious generalization of \eqref{constraints_inter}. 
%In addition, we assume $\bar{u}_{12}=\bar{u}_{21} = 0$ and $\bar{T}_{12}= \bar{T}_{21}$ meaning that
%\begin{align}
%\begin{split}
%\frac{1}{\bar{n}_1} \int \nu_{12} v M_{\lambda_{12}} dv = \frac{1}{\bar{n}_2} \int \nu_{21} v M_{\lambda_{21}} dv\\
%\frac{3 m_1}{\bar{n}_1} \int \nu_{12} |v|^2 M_{\lambda_{12}} dv = \frac{3 m_2}{\bar{n}_2} \int \nu_{21} |v|^2 M_{\lambda_{21}} dv
%\end{split}
%\label{equal_vel_temp}
%\end{align}
%As in the two species case, the form \eqref{oneN} can be motivated by considering the minimization problem
%\begin{align}
%\min_{g \in \chi_i} \int \nu_{ii} h(g) dv, ~i=1,...,N
%\end{align}
%with 
%\begin{align}
%\chi_i= \lbrace g | g>0, \nu_{ii} (1+ |v|^2) g \in L^1(\mathbb{R}^3),\int \nu_{ii} a^i(v) (g_i - f_i) dv = 0  \rbrace
%\end{align}
%for the one species relaxation terms. The mixture Maxwellians \eqref{mixN} are defined by considering the minimization problem
%\begin{align}
%\min_{g_{ij}, g_{ji} \in \chi}  \int \nu_{ij} h(g_{ij}) dv +  \int \nu_{ji} h(g_{ji}) dv , 
%\end{align}
%with 
%\begin{align}
%\begin{split}
%\chi= \lbrace g_{ij}, g_{ji} | g_{ij}, g_{ji} >0, \nu_{ij} (1+ |v|^2) g_{ij}, \nu_{ji} (1+ |v|^2) g_{ji}  \in L^1(\mathbb{R}^3),\\ \int \nu_{ij} g_{ij} dv = \int \nu_{ij} f_i dv, \quad \int \nu_{ji} g_{ji} dv = \int \nu_{ji} f_j dv \\
%\int \nu_{ij} \begin{pmatrix}
%m_i v \\ \frac{m_i}{2} |v|^2
%\end{pmatrix} (g_{ij} - f_i) dv  = - 
%\int \nu_{ji} \begin{pmatrix}
%m_j v \\ \frac{m_i}{2} |v|^2
%\end{pmatrix} (g_{ji} - f_j) dv \rbrace
%\end{split}
%\end{align}
%for each pair $(i,j)$. 
All the proofs concerning existence and uniqueness of the target Maxwellians and the H-Theorem can be proven exactly in the same way as for two species. %For the H-Theorem, we prove
%\begin{align}
%\int \ln f_i Q_{ij}(f_i.f_j) dv + \int \ln f_j Q_{ji}(f_j, f_i) dv \leq 0
%\label{HN}
%\end{align}
%for each pair $(i,j)$ in the same way as for $(12)$ in the two species case. 
For the total entropy $H(f_1,...,f_N) = \int (h(f_1)  + \dots + h(f_N)) dv$ we obtain
\begin{align}
\partial_t \left(  H(f_1,..., f_N) \right) + \nabla_x \cdot \left(\int  v (h(f_1) + \dots + h(f_N)) dv \right) \leq 0.
\end{align}
%One sees this by multiplying the equation for $f_j, ~ j=1,...,N$ by $\ln f_j$ and using \eqref{HN} for all species pairs.

\section*{Conclusion}
We have presented a multi-species BGK model in which the collision frequencies depend on the microscopic velocity.  The model is formally derived based on an entropy minimization principle, which implies that the target functions take the form of Maxwellians. However, contrary to classical BGK models with velocity-independent frequencies, the relationship between the Maxwellian parameters and the moments of the distribution function is not analytic.
Thus some effort is required to establish rigorously the existence of parameters which satisfy first-order optimality conditions.  We also show that the derived model satisfies an H-Theorem and that it can be extended to the case of arbitrarily many species undergoing binary collisions.

In future work, we will develop numerical tools for discretizing the model developed here, including the numerical solution of the defining optimization problem.  A numerical code will enable computational explorations about how to choose the collision frequencies and what benefit is providing by their flexibility.  Also, because the motivation for the model is the simulation of multi-species plasmas, we will extend it for use in such contexts by adding self-consistent fields.

\section*{Acknowledgements}
Christian Klingenberg acknowledges a grant by the Bayrische Forschungsallianz. \\ \\
Marlies Pirner is supported from the Humboldt foundation and from the Austrian Science Fund (FWF) through grant number F65. \\ \\
The work of Jeff Haack was supported by the US Department of Energy through the Los Alamos National Laboratory. Los Alamos National Laboratory is operated by Triad National Security, LLC, for the National Nuclear Security Administration of U.S. Department of Energy (Contract No. 89233218CNA000001). Los Alamos Report LA-UR-20-21464. \\ \\
The work of Cory Hauck is sponsored by the Office of Advanced Scientific Computing Research, U.S. Department of Energy, and performed at the Oak Ridge National Laboratory, which is managed by UT-Battelle, LLC under Contract No. De-AC05-00OR22725 with the U.S. Department of Energy. The United States Government retains and the publisher, by accepting the article for publication, acknowledges that the United States Government retains a non-exclusive, paid-up, irrevocable, world-wide license to publish or reproduce the published form of this manuscript, or allow others to do so, for United States Government purposes. The Department of Energy will provide public access to these results of federally sponsored research in accordance with the DOE Public Access Plan (http://energy.gov/downloads/doe-public-access-plan).

% \bibitem{Pirner} C. Klingenberg, M.Pirner, G.Puppo, \emph{A consistent kinetic model for a two-component mixture with an application to plasma}, Kinetic and related Models, 2017 

%\end{thebibliography}
%\textcolor{magenta}{Regularity of the limit of a smooth sequence converging in $L^1$?, proof only possible for classical solutions, otherwise if function is only in $L^1$ I can not fix a certain point}
%\newpage

\end{document}